\documentclass[reqno,11pt]{amsart} 
\usepackage{amsmath}
\usepackage{amssymb}
\usepackage{amsthm}
\usepackage{verbatim}
\usepackage{xcolor}
\usepackage{graphics}

\newcommand\NoBlackBoxes{\global\overfullrule0pt}

\NoBlackBoxes
\theoremstyle{plain}

\begin{document}

\title{A QUANTITATIVE CRAM\'ER-WOLD THEOREM \\
FOR ZOLOTAREV DISTANCES
}

\author{Sergey G. Bobkov$^{1,3}$}
\thanks{1) 
School of Mathematics, University of Minnesota, Minneapolis, MN, USA,
bobkov@math.umn.edu. 
}

\author{Friedrich G\"otze$^{2,3}$}
\thanks{2) Faculty of Mathematics,
Bielefeld University, Germany,
goetze@math-uni.bielefeld.de.
}

\thanks{3) Research supported by the NSF grant DMS-2154001 and
the GRF – SFB 1283/2 2021–317210226}


\begin{abstract}
An upper bound for Zolotarev distances between probability
measures on multidimensional Euclidean spaces is given in terms of similar
distances between one dimensional projections.
\end{abstract}

\maketitle
\markboth{Sergey G. Bobkov and Friedrich G\"otze}{Cram\'er-Wold theorem}

\def\theequation{\thesection.\arabic{equation}}
\def\E{{\mathbb E}}
\def\R{{\mathbb R}}
\def\C{{\mathbb C}}
\def\P{{\mathbb P}}
\def\Z{{\mathbb Z}}
\def\L{{\mathbb L}}
\def\T{{\mathbb T}}

\def\G{\Gamma}

\def\Ent{{\rm Ent}}
\def\var{{\rm Var}}
\def\Var{{\rm Var}}

\def\H{{\rm H}}
\def\Im{{\rm Im}}
\def\Tr{{\rm Tr}}
\def\s{{\mathfrak s}}

\def\k{{\kappa}}
\def\M{{\cal M}}
\def\Var{{\rm Var}}
\def\Ent{{\rm Ent}}
\def\O{{\rm Osc}_\mu}

\def\ep{\varepsilon}
\def\phi{\varphi}
\def\vp{\varphi}
\def\F{{\cal F}}

\def\be{\begin{equation}}
\def\en{\end{equation}}
\def\bee{\begin{eqnarray*}}
\def\ene{\end{eqnarray*}}

\thispagestyle{empty}

\section{{\bf Introduction}}
\setcounter{equation}{0}

\vskip2mm
\noindent
A sequence of probability measures $(\mu_n)_{n \geq 1}$ is weakly convergent 
to a probability measure $\mu$ on $\R^d$, writing for short $\mu_n \Rightarrow \mu$, if
$$
\int_{\R^d} u\,d\mu_n \rightarrow \int_{\R^d} u\,d\mu \quad {\rm as} \ n \rightarrow \infty
$$
for any bounded continuous function $u$ on $\R^d$. According to the
Cram\'er-Wold continuity theorem, this property is equivalent to the assertion that 
$\mu_{n,\theta} \Rightarrow \mu_\theta$ on the real line for all one dimensional projections, 
that is, for the images $\mu_{n,\theta}$ and $\mu_\theta$ of these measures under linear 
mappings $x \rightarrow \left<x,\theta\right>$ from $\R^d$ to $\R$ with an arbitrary fixed
$\theta \in \R^d$ (cf. \cite{C-W}, \cite{Bil}).

The  problem of how to quantify this characterization by means of various distances
describing  weak convergence is interesting and not clear. Another natural problem is 
to determine whether or not there is a similar description for stronger metrics, which 
could potentially reduce a number of high dimensional questions
to dimension one, perhaps under proper moment assumptions. 
In connection with computer tomography problems and continuity properties 
of the Radon transform (which assigns to every probability distribution $\mu$ on $\R^d$ 
the collection of all its one-dimensional projections $\mu_\theta$), this question and related inversion issues 
were discussed in 1990's in the works by Zinger, Klebanov and Khalfin (cf. e.g. \cite{K-K}).
A natural choice in this problem is given by the Kolmogorov distance
$$
\rho_d(\mu,\nu) = \sup |\mu(H) - \nu(H)|,
$$
where the supremum is running over all half-spaces $H$ in $\R^d$. 
However, it was shown by Zaitsev \cite{Z} that the smallness of 
$\sup_{|\theta| = 1} \rho_1(\mu_\theta,\nu_\theta)$ does not guarantee that 
$\rho_d(\mu,\nu)$ with $d \geq 2$ will be small as well, even if these measures are compactly 
supported. Moreover, here $\rho_1$ may be strengthened to the total variation distance.
The corresponding counter-example shows that $\rho_d$ is essentially stronger
than any metric which metrizes the topology of weak convergence in the space
of probability distributions on $\R^d$.

With respect to the Kantorovich transport distance
$$
W(\mu,\nu) = \inf \int_{\R^d}\int_{\R^d} |x-y|\,d\pi(x,y)
$$
the problem has been recently considered in \cite{B-G}. Here the infimum is taken 
over all probability measures $\pi$ on $\R^d \times \R^d$ with marginals $\mu$ and $\nu$.  
Denote by ${\mathfrak P}_q$ and ${\mathfrak P}_q(b)$, $q \geq 1$, the collections of all
probability measures $\mu$ on $\R^d$ such that $\int |x|^q\,d\mu(x) < \infty$
and respectively 
$$
\int |x|^q\,d\mu(x) \leq b^q,
$$ 
where $b > 0$ is a parameter and $|x|$ denotes the Euclidean norm of $x$.
By a well-known characterization, $W$ metrizes the topology of weak convergence in the 
space ${\mathfrak P}_1$ (\cite{V}, Theorem 7.12). Regarding the Cram\'er-Wold theorem,
it was shown in \cite{B-G} that
restricting $\mu$ and $\nu$  to the class ${\mathfrak P}_q$ with $q>1$, the smallness of 
$\sup_{|\theta| = 1} W(\mu_\theta,\nu_\theta)$ does guarantee that $W(\mu,\nu)$ is small. 
More precisely, we have the following quantified version of this theorem.

\vskip3mm
{\bf Theorem 1.1.} {\sl For all probability measures $\mu$ and $\nu$ in ${\mathfrak P}_q(b)$
with $q > 1$ and $b \geq 0$,
\be
W(\mu,\nu) \leq c\,b^{1-\beta} \sup_{|\theta| = 1}  W(\mu_\theta,\nu_\theta)^\beta, \quad \
\beta = \frac{2}{2 + d q/(q-1)},
\en
where $c>0$ is an absolute constant.
}

\vskip5mm
For compactly supported measures, for example, when both $\mu$ and $\nu$ are supported on the
Euclidean ball $|x| \leq 1$, this bound can be simplified to
\be
W(\mu,\nu) \leq c\,\sup_{|\theta| = 1}  W(\mu_\theta,\nu_\theta)^\beta, \quad \ \beta = \frac{2}{2 + d}.
\en

In more recent studies, the supremum in (1.1) is often called the max-scliced Wasserstein 
distance (but more correctly it should be called the max-scliced Kantorovich distance). 
Restricted to the class of empirical measures on $\R^d$, this quantity reflects the uniform 
behaviour of their one-dimensional projections, which represent empirical measures
on the real line. It has been the subject of many recent investigations towards the problem 
of convergence rates, cf. e.g.  \cite{NW-R}, \cite{N-G-S-K}, \cite{M-B-W}, \cite{Boe} 
(including the Kantorovich distances $W_p$ generated by the power cost functions $|x-y|^p$).
In particular, an application of (1.2) to empirical measures shows that necessarily $\beta \leq 2/d$
for the upper bound in (1.2). Hence the present power $\beta = 2/(2 + d)$ 
is asymptotically optimal with respect to the growing dimension $d$.

In this note we consider the class of Zolotarev ideal metrics which include  $W = W_1$ 
as a partial case. Given two measures $\mu$ and $\nu$ in ${\mathfrak P}_p$, 
the Zolotarev distance of an integer order $p \geq 1$ between $\mu$ and $\nu$ is defined as
\be
\zeta_p(\mu,\nu) \, = \,  \sup \Big|\int_{\R^d} u\,d\mu - \int_{\R^d} u\,d\nu \Big|
\en
where the supremum is taken over all $p$ times differentiable (or equivalently, $C^\infty$-smooth) 
functions $u:\R^d \to \R$ whose partial derivatives of order $p$ are bounded by 1 
in absolute value. These distances were introduced and explored by Zolotarev in a series 
of papers in the 1970’s, especially in the context of the central limit theorem, 
cf. \cite{Z1}, \cite{Z2}, \cite{Z3}.

If $p=1$, this quantity is equivalent to $W$, by the Kantorovich duality theorem.
In the case $d=1$, in contrast to the multidimensional
situation, it has a simple description
$$
\zeta_1(\mu,\nu) = W(\mu,\nu) = \int_{-\infty}^\infty |F_\mu(x) - F_\nu(x)|\,dx,
$$
where $F_\mu(x) = \mu((-\infty,x])$, $F_\nu(x) = \nu((-\infty,x])$ are associated distribution 
functions. However, when $p>1$, although the integrals in (1.3) are well-defined and
finite, the finiteness of the supremum requires that the measures $\mu$ and $\nu$ have 
to have  equal mixed moments up to order $p-1$. For example,
the most frequently used distance $\zeta_3$ requires that $\mu$ and $\nu$
have equal means (baricenters) and equal covariance matrices.

Since the functions $u(x) = v(\left<x,\theta\right>)$ with $|\theta|=1$ and
$v:\R \rightarrow \R$ such that $|v^{(p)}(t) \leq 1$ for all $t \in \R$
participate in the supremum (1.3), we necessarily get
$$
\zeta_p(\mu,\nu) \geq \sup_{|\theta| = 1} \zeta_p(\mu_\theta,\nu_\theta).
$$
Here we reverse this inequality in a somewhat similar form
in analogy with Theorem 1.1.

\vskip5mm
{\bf Theorem 1.2.} {\sl If the probability measures $\mu$ and $\nu$ belong to 
${\mathfrak P}_q(b)$ with $q > p$ and $b \geq 0$, then with some absolute constant $c>0$ 
\be
\zeta_p(\mu,\nu) \leq (cd)^p\,b^{1-\beta} \sup_{|\theta| = 1}  \zeta_p(\mu_\theta,\nu_\theta)^\beta,
\en
where
$$
\beta = \frac{2}{2 + dq/p(q-p)}.
$$
}

\vskip2mm
This inequality makes sense when $\mu$ and $\nu$ have equal mixed moments 
up to order $p-1$, in which case both sides of (1.4) are finite. 
Letting $q \rightarrow \infty$ for compactly supported measures, the relation (1.4) 
can be simplified as follows.

\vskip5mm
{\bf Corollary 1.3.} {\sl If $\mu$ and $\nu$ are supported on the unit ball $|x| \leq 1$ in $\R^d$,
then
$$
\zeta_p(\mu,\nu) \leq (cd)^p\, \sup_{|\theta| = 1}  \zeta_p(\mu_\theta,\nu_\theta)^\beta \quad
{\sl with} \ \ \beta = \frac{2}{2 + d/p}.
$$
}

\vskip2mm
Similarly to the proof Theorem 1.1, the proof of Theorem 1.2 involves truncation 
and Fourier analysis. However, there are additional obstacles  which require the construction 
of special signed measures for a  smoothing procedure and the study  of their basic properties. 
Thus the plan of this note is as follows.

\vskip2mm
{\sl Contents}:

\vskip2mm
1. Introduction.

2. Zolotarev semi-norms.

3. Smoothing.

4. Special compactly supported functions.

5. Reduction to compactly supported functions.

6. Fourier transforms.

7. Proof of Theorem 1.2.

\vskip7mm
\section{{\bf Zolotarev Semi-Norms}}
\setcounter{equation}{0}

\vskip2mm
\noindent
Since the smoothing will involve signed measures, it is desirable to extend the definition (1.3)
to this larger class. 
Let $\lambda$ be a signed (Borel) measure on $\R^d$ with total mass $\lambda(\R^d) = 0$
such that its variation $|\lambda|$ (as a positive measure) satisfies $\int |x^p|\,d|\lambda|(x) < \infty$
(the latter property is referred in the sequel as finiteness of the $p$-th absolute moment). Put
\be
\|\lambda\|_{\zeta_p} \, = \,  \sup \Big|\int_{\R^d} u\,d\lambda\Big|,
\en
where the supremum is taken over all $p$ times differentiable (or equivalently, $C^\infty$-smooth) 
functions $u:\R^d \to \R$ whose partial derivatives of order $p$ are bounded by 1.
This quantity may be called Zolotarev semi-norm of $\lambda$ of an integer order $p \geq 1$.
Thus,
$$
\zeta_p(\mu,\nu) = \|\mu - \nu\|_{\zeta_p}.
$$

In the sequel we denote by $\alpha = (p_1,\ldots,p_d)$ a multi-index with $p_i \geq 0$ 
being integers and by  $|\alpha| = p_1 + \cdots + p_d$ its length.
Let us write $\alpha! = p_1! \dots p_d!$\, and
$$
x^\alpha = x_1^{p_1} \dots x_1^{p_d} \quad {\rm for} \ x = (x_1,\dots,x_d) \in \R^d
$$
and furthermore
$$
D^\alpha u \, = \, 
\frac{\partial^{|\alpha|}}{\partial x_1^{p_1} \cdots \partial x_d^{p_d}}\, u 
$$
for the partial derivatives of $u$ of order $|\alpha|$, 
using the convention that $D^\alpha u = u$ when $\alpha = (0,\dots,0)$.

It is easily seen  (cf. Corollary 2.4 for one direction) that
$\|\lambda\|_{\zeta_p}$ is finite if and only if the measure $\lambda$ has 
finite $p$-th absolute moment and zero mixed moments up to order $p-1$, that is,
\be
\int_{\R^d} x^\alpha\,d\lambda(x) = 0, \quad |\alpha| \leq p-1.
\en
When $\alpha = (0,\dots,0)$, this is the same as $\lambda(\R^d) = 0$,
which is the only condition we require for $p=1$.

\vskip5mm
{\bf Lemma 2.1.} {\sl If the signed measure $\lambda$ on $\R^d$ with finite $p$-th absolute
moment satisfies the moment condition $(2.2)$, then the supremum in $(2.1)$ 
may be restricted to $C^\infty$-smooth functions $u$ on $\R^d$ such that
\be
D^\alpha u(0) = 0 \quad {\sl for \ all} \ \alpha\  {\sl with} \ |\alpha| \leq p-1
\en 
and
\be
\sup_x\, |D^\alpha u(x)| \leq 1 \quad {\sl for \ all} \ \alpha\  {\sl with} \ |\alpha| = p.
\en
}

\vskip3mm
Indeed, if $u$ satisfies (2.3)-(2.4), so does
$$
\widetilde u(x) = u(x) - \sum_{|\alpha| \leq p-1} \frac{D^\alpha u(0)}{\alpha!}\,x^\alpha.
$$
Due to (2.2), we also have $\int \widetilde u\,d\lambda = \int u\,d\lambda$, thus proving
the lemma.

We will need an elementary non-uniform bound on partial derivatives of the
functions $u$ participating in the supremum (2.1) subject to the constraints of Lemma 2.1.

\vskip5mm
{\bf Lemma 2.2.} {\sl Given a $C^p$-smooth function $u$ on $\R^d$  satisfying $(2.3)-(2.4)$,
we have, for any multi-index $\alpha$ of length $m \leq p$,
\be
|D^\alpha u(x)| \leq \frac{1}{(p - m)!}\,(|x_1| + \dots + |x_d|)^{p-m}, \quad x \in \R^d.
\en
}

\vskip2mm
{\bf Proof.} The case $m=p$ reduces (2.5) to (2.4). If $m \leq p-1$, then according to 
the Taylor multi-dimensional integral formula, applied to the $\alpha$-th 
partial derivative of $u$ around the origin, we have, for any point $x = (x_1,\dots,x_d) \in \R^d$, 
\bee
D^\alpha u(x) 
 & = & 
\sum_{|\beta| \leq p-m-1} \frac{D^{\alpha + \beta} u(0)}{\beta!}\,x^\beta \\
 & & + \
\sum_{|\beta| = p-m} \frac{|\beta|}{\beta!}\, x^\beta \int_0^1 (1-t)^{p-m-1} D^{\alpha + \beta} u(tx)\,dt,
\ene
where the summation extends over all multi-indices $\beta$ with
lengths as indicated in the summation. The first sum is vanishing, by (2.3). By (2.4), we also have 
$|D^{\alpha + \beta} u(tx)| \leq 1$ in the second sum, implying that the last integral does not
exceed $1/|\beta|$. Therefore, applying the multinomial formula, we get
$$
|D^\alpha u(x)| \leq \sum_{|\beta| = p-m} \frac{1}{\beta!}\, |x^\beta| = \frac{1}{(p-m)!}\,
(|x_1| + \dots + |x_d|)^{p-m}.
$$
\qed

\vskip3mm
{\bf Remark 2.3.} When dropping the condition $u(0)=0$, the bound (2.5) continues to hold
for partial derivatives of orders $1 \leq m \leq p$. In this case the condition (2.3)
needs to hold for $1 \leq |\alpha| \leq p-1$ only.

\vskip5mm
{\bf Corollary 2.4.} {\sl If the signed measure $\lambda$ on $\R^d$ with finite $p$-th absolute
moment satisfies the moment condition $(2.2)$, then
$$
\|\lambda\|_{\zeta_p} \, \leq \, \frac{1}{p!}\,d^{\frac{p}{2}} \int |x|^p\,d|\lambda|(x).
$$
}

\vskip2mm
Indeed, applying Lemma 2.2 with $m=0$, we obtain that
\be
|u(x)| \leq \frac{1}{p!}\,(|x_1| + \dots + |x_d|)^p \leq \frac{1}{p!}\,d^{\frac{p}{2}}\,|x|^p, 
\quad x \in \R^d.
\en
It remains to integrate this inequality over the measure $|\lambda|$, to recall the definition (2.1)
and to apply Lemma 2.1.

\vskip7mm
\section{{\bf Smoothing}}
\setcounter{equation}{0}

\vskip2mm
\noindent
The smoothing by means of convolution can considerably improve analytic properties of measures,
which may be expressed in terms of the associated Fourier transforms.
Forced by the requirement of vanishing mixed moments, we  need convolution
by signed measures $\kappa$ which are supported
on the unit ball $B_1 = \{x \in \R^d: |x| \leq 1\}$ such that
\be
\kappa(B_1) = 1, \quad
\int_{B_1} x^\alpha\,d\kappa(x) = 0 \quad {\rm for} \ 1 \leq |\alpha| \leq p-1,
\en 
where $\alpha$ denote multi-indices. Denote by $\kappa_\ep$ the image 
of $\kappa$ under the linear map $x \rightarrow \ep x$ with parameter $\ep>0$.
From (3.1) it follows that, for any signed measure $\lambda$ on $\R^d$ with finite
$p$-th absolute moment, the convolution 
\be
\lambda_\ep = \lambda * \kappa_\ep
\en
has the same mixed moments as $\lambda$ up to order $p-1$.

\vskip5mm
{\bf Proposition 3.1.} {\sl Given a signed measure $\lambda$ on $\R^d$ with finite
$p$-th absolute moment, for any $\ep>0$, 
\be 
\|\lambda_\ep - \lambda\|_{\zeta_p} \, \leq \, 
\frac{(d\ep)^p}{p!}\, \|\lambda\|_{\rm TV} \, \|\kappa\|_{\rm TV}. 
\en
}

\vskip2mm
Here and in the sequel, we denote by $\|\lambda\|_{\rm TV}$ the total variation norm
of $\lambda$.

\vskip5mm
{\bf Proof.}
Let $u$ be a $C^\infty$-smooth function on $\R^d$ whose partial derivatives of order $p$
are bounded by 1 in absolute value. By definition of the convolution, and using the assumption 
$\kappa(B_1) = 1$, we have
\be
\int_{\R^d} u\,d \lambda_\ep - \int_{\R^d} u\,d\lambda = 
\int_{\R^d} \int_{B_1} \big (u(x+\ep y) - u(x) \big)\,d\lambda(x)\,d\kappa(y).
\en
In order to bound the last integral, we apply a multidimensional integral Taylor formula at the point $x$
and write
$$
u(x+\ep y) - u(x) \, = \, 
\sum_{1 \leq |\alpha| \leq p-1} \frac{D^\alpha u(x)}{\alpha!}\, \ep^{|\alpha|}\,y^\alpha +
\ep^p \sum_{|\alpha| = p} Q_\alpha\, y^\alpha,
$$
where
$$
Q_\alpha \, = \,  \frac{|\alpha|}{\alpha!} \int_0^1 (1-t)^{p - 1}\, 
D^\alpha u(x+t \ep y)\,dt.
$$
Due to the moment assumption in (3.1), the integration over $\kappa$ yields 
$$
\int_{B_1} (u (x+\ep y)  - u(x))\,d\kappa(y)  \, = \, \ep^p
\sum_{|\alpha| = p} \frac{|\alpha|}{\alpha!} \int_0^1 (1-t)^{p-1} 
\bigg[\int_{B_1} D^\alpha u(x+t \ep y)\,y^\alpha\,d\kappa(y)\bigg] dt. 
$$
Here $|y^\alpha| \leq 1$ for $y \in B_1$, while $|D^\alpha u(x+t \ep y)| \leq 1$.
Hence the last integral does not exceed $\|\kappa\|_{\rm TV}$ in absolute value, leading to
$$
\Big|\int_{B_1} (u(x+\ep y) - u(x))\,d\kappa(y)\Big| \, \leq \, \ep^p\, \|\kappa\|_{\rm TV}
\sum_{|\alpha| = p}  \frac{1}{\alpha!} \, = \,  \ep^p \|\kappa\|_{\rm TV} \cdot \frac{d^p}{p!}.
$$
Using this bound in (3.4), it remains to take the supremum over all admissible functions $u$,
according to the definition (2.1). 
\qed

\vskip5mm
In view of (3.1), it is natural to choose $\kappa$ with total variation norm as small as possible 
and  constructing some concrete 
examples. We will also need to control the maximum of the density of $\kappa$.

\vskip5mm
{\bf Lemma 3.2.} {\sl For any integer $p \geq 1$, there exists a rotationally invariant signed 
measure $\kappa$ on $B_1$, which satisfies $(3.1)$ and has density $w$ such that
\be
|w(x)| \leq \frac{1}{\omega_d}\,p^{d+1}\,8^p, \quad x \in B_1.
\en
A similar bound holds for the total variation norm:  $\|\kappa\|_{\rm TV} \leq p\,2^{3p-1}$.
}

\vskip5mm
Here and elsewhere $\omega_d$ stands for the $d$-dimensional volume of the unit ball $B_1$. 

Applying Lemma 3.2 in Proposition 3.1, we get:

\vskip5mm
{\bf Corollary 3.3.} {\sl Given a signed measure $\lambda$ on $\R^d$ with finite
$p$-th absolute moment, we have, for any $\ep>0$, 
\be 
\|\lambda_\ep - \lambda\|_{\zeta_p} \, \leq \, \frac{(8d\ep)^p}{2(p-1)!}\, \|\lambda\|_{\rm TV},
\en
where the convolution $\lambda_\ep$ is defined in $(3.2)$ for the measure $\kappa$ from Lemma $3.2$.
}

\vskip5mm
In dimension $d=1$, the measures $\kappa$ as in Lemma 3.2 together with integral bounds on their
Fourier transforms were studied in \cite{B}, and recently in \cite{B-L} for the class of
product measures on the torus.
For completeness, we borrow basic steps in the construction of $\kappa$.

\vskip5mm
{\bf Proof of Lemma 3.2.} If $p=1$, one may choose the uniform distribution $\kappa$ on $B_1$,
that is, with density $w(x) = \frac{1}{\omega_d}\,1_{B_1}(x)$. So, let $p \geq 2$. First we construct 
a signed measure $\kappa_0$ supported on the interval $(0,1)$ such that
\be
\kappa_0(0,1) = 1, \quad \int_0^1 r^k\,d\kappa_0(r) = 0, \quad k = 1,\dots,p-1.
\en
Identifying distribution functions with measures, take a probability measure $H$ on $(0,1)$ and 
define
$$
\kappa_0(r) \,= \, a_1 H\Big(\frac{r}{b_1}\Big) + \cdots + a_p H\Big(\frac{r}{b_p}\Big)
$$
with weights $a_1, \ldots, a_p \in \R$ and parameters $0 < b_1 < \cdots < b_p \leq 1$. 
Then $\kappa_0$ is supported on $(0,1)$ as a measure and has total variation norm
\be
\|\kappa_0\|_{\rm TV} \,\leq  \, \sum_{i=1}^p |a_i|.
\en
For such a measure $\kappa_0$, the condition (3.7) is fulfilled as long as
$$
\sum_{i=1}^p a_i = 1, \quad \sum_{i=1}^p a_i b_i^k  =  0, \quad k = 1, \ldots, p-1.
$$
This is a linear system of $p$ equations in $p$ unknowns $a = (a_1,\ldots,a_p)$, 
which can be written in matrix form as $Va = e_1$, where $V$ is the Vandermonde matrix
\begin{displaymath}
V \ = \ 
\left(\begin{array}{cccc}
1          & 1         & \cdots  & 1         \\
b_1       & b_2     & \cdots  & b_p      \\
\vdots  & \vdots  & \ddots & \vdots  \\
b_1^{p-1}  & b_2^{p-1}  & \cdots  & b_p^{p-1} \\
\end{array}
\right)
\end{displaymath}
and $e_1= (1,0,\ldots,0)$ (as a column). It has a non-zero determinant
$ \mathrm {det}(V) = \prod_{i<j}\, (b_j - b_i)$,
so that $V$ is invertible and $a = V^{-1} e_1$. 

In order to estimate the coordinates of the vector $a$ in terms of $b_i$'s,
one may employ the following result from \cite{G} on the inverse of the Vandermonde 
matrix. Let us define the norm of a $p \times p$ matrix $M = (m_{ij})$ by
$$
\|M\| = \max_{1 \leq i \leq p}\, \sum_{j=1}^p |m_{ij}|.
$$
Then we have
$$ 
\|V^{-1}\| \, \leq \, \max_{1 \leq i \leq p}  \, \prod_{j \neq i} \frac{1 + |b_j|}{|b_i - b_j|} \, .
$$
In particular, the choice $b_i = i/p$ which we use below leads to
\bee
\|V^{-1}\| 
  & \leq & 
\max_{1 \leq i \leq p}\, 
\prod_{j=1, \ j \neq i}^p\,\frac{p + j}{|j - i|}
  \, = \,
\max_{1 \leq i \leq p}\, \frac{(p+1)\cdots 2p}{(p+i)\,(i-1)!\,(p-i)!} \\
  & = & 
2p\, \binom{2p-1}{p} \max_{1 \leq i \leq p}\, \frac{\binom {p-1}{i-1}}{p+i} \, \leq \, 2^{3p-1},
\ene
where we made use of the bound $\binom {n}{i} \leq 2^n$ for the binomial coefficients. 

Now, since $a_i = (V^{-1} e_1)_i = (V^{-1})_{i1}$, it follows that 
$$
|a_i| \leq \|V^{-1}\| \leq 2^{3p-1},
$$ 
and by (3.8), 
\be
\|\kappa_0\|_{\rm TV} \,\leq  \, p\,2^{3p-1}.
\en

This bound does not depend on the choice of $H$. Let us choose the measure
with distribution function $H(r) = r^d$, $0 \leq r \leq 1$, which defines 
$\kappa_0$ with distribution function
$$
\kappa_0(r) = \sum_{i=1}^p a_i \Big(\frac{r}{b_i}\Big)^d\,1_{[0,b_i]}(r),
\quad 0 \leq r \leq 1.
$$
Then, the desired measure $\kappa$ may be defined as the image of the product 
measure $\kappa_0 \otimes \sigma_{d-1}$ under the mapping 
$(r,\theta) \rightarrow r \theta$, where $\sigma_{d-1}$ is the uniform distribution
on the unit sphere $S^{d-1}$ in $\R^d$ (being the Bernoulli measure 
$\frac{1}{2}\,\delta_1 + \frac{1}{2}\,\delta_{-1}$ for dimension $d=1$). In particular, 
$\kappa(B_1) = 1$, and, by (3.9),
$$
\|\kappa\|_{\rm TV} \, \leq \, \|\kappa_0 \otimes \sigma_{d-1}\|_{\rm TV}
\, \leq \, \|\kappa_0\|_{\rm TV} \,\leq  \, p\,2^{3p-1}.
$$
Thus, we obtain the second claim of the lemma about the total variation norm.

To prove the first claim, note that, for every $i = 1,\dots,p$, the same mapping pushes forward
the product probability measure $H(r/b_i) \otimes \sigma_{d-1}$ with the marginal
distribution function $H(r/b_i) = r^d/b_i^d$ ($0 \leq r \leq b_i$) to the uniform distribution
on the ball $|x| < b_i$, with density
$$
w_i(x) = \frac{1}{\omega_d b_i^d}\,1_{\{|x| < b_i\}}, \quad x \in \R^d.
$$ 
Hence $\kappa$ has density
$$
w(x) \, = \, \sum_{i=1}^p a_i w_i(x) \, = \,
\frac{p^d}{\omega_d}\, \sum_{i=1}^p \frac{a_i}{i^d}\,1_{\{|x| < i/p\}}, 
$$
implying that
$$
|w(x)| \, \leq \,
\frac{p^d}{\omega_d}\,\sum_{i=1}^p \frac{1}{i^d}  \max_i |a_i| \, \leq \, 
\frac{2p^d}{\omega_d}\,2^{3p-1}, \quad d \geq 2.
$$
In the case $d=1$, we similarly have
$$
|w(x)| \, \leq \,
\frac{p}{\omega_1}\,\sum_{i=1}^p \frac{1}{i}  \max_i |a_i| \, \leq \, 
\frac{p}{\omega_1}\,2^{3p-1} (1 + \log p) \, \leq \, \frac{p^2}{\omega_1}\,2^{3p-1}.
$$
This proves (3.5). Finally, 
\bee
\int_{B_1} x^\alpha\,d\kappa(x) 
 & = &
\int_0^1 \!\!\int_{S^{d-1}} (r \theta)^\alpha\,d\kappa_0(r)\,d\sigma_{d-1}(\theta) \\
 & = &
\int_0^1 r^{|\alpha|}\,d\kappa_0(r) 
\int_{S^{d-1}} \theta^\alpha\,d\sigma_{d-1}(\theta) \, = \, 0
\ene
for any multi-index $\alpha$ such that $1 \leq |\alpha| \leq p-1$, by the property (3.7).
Thus, the requirement (3.1) is fulfilled.
\qed

\vskip7mm
\section{{\bf Special Compactly Supported Functions}}
\setcounter{equation}{0}

\vskip2mm
\noindent
The next step is to prepare the restriction of the supremum in (2.1) to the class of 
compactly supported functions. It is based on some special functions of independent interest
as stated in the following:

\vskip5mm
{\bf Lemma 4.1.} {\sl There exists a $C^p$-smooth function $\psi:\R^d \rightarrow [0,1]$ such that

\vskip2mm
$1)$ \ $\psi(x) = 1$ for $|x| \leq \frac{1}{2}$;

$2)$ \ $\psi(x) = 0$ for $|x| \geq 1$;

$3)$ \ For any multi-index $\gamma$
of length $m \leq p$,
\be
\sup_x |D^\gamma \psi(x)| \leq c^p p^m
\en
for some absolute constant $c$.
}

\vskip5mm
The constant on the-right-hand of (4.1) does not depend on $d$. This  is reflected in 
the dimension dependence of the constant used in the bound of Theorem 1.2.

\vskip5mm
{\bf Proof.}
The  function $\psi$  may be chosen as $\psi(x) = \xi(|x|^2)$ with
$\xi:[0,\infty) \rightarrow [0,1]$ being  $C^p$-smooth and satisfying $\xi(t) = 1$ for $t \leq \frac{1}{4}$
and $\xi(t) = 0$ for $t \geq 1$ in order to meet the requirements 1)-2). More precisely, define
\be
v(t) = \frac{1}{A_p} \int_0^t (s(1-s))^p\,ds, \quad 0 \leq t \leq 1, \qquad
v(t) = 1 \quad t \geq 1,
\en
where $A_p$ is a normalizing constant such that $v(1)=1$.
All derivatives of $v(t)$ at $t=0$ and $t=1$ are vanishing up to order $p$, which guarantees
the $C^p$-smoothness of this function on the positive half-axis $t \geq 0$. Then put
$\xi(t) = 1 - v(\frac{4t-1}{3})$, $t \geq \frac{1}{4}$ and  extend it via  $\xi(t) = 1$ to the interval 
$0 \leq t \leq \frac{1}{4}$. Thus, we arrive at a function $\psi$ on $\R^d$ given by
\be
\psi(x) = 1 -  v(y), \quad y = \frac{4\,|x|^2 - 1}{3}, \ \ x = (x_1,\dots,x_d),
\en
and need to give a uniform bound (4.1) in the region $\frac{1}{2} < |x| < 1$.
If $m = 0$, then, by definition, $0 \leq v(t) \leq 1$, hence $0 \leq \psi(x) \leq 1$ on the unit ball $|x| \leq 1$.

In the case $m \geq 1$, the differentiation can be based on the multinomial extension of
the Newton formula to the product of two or more functions:
If $f = f_1 \dots f_n$ on $\R^d$, then
\be
D^\gamma f = \sum 
\binom{|\gamma_1|}{\gamma_1} \dots  \binom{|\gamma_n|}{\gamma_n} 
D^{\gamma_1} f_1 \dots D^{\gamma_n} f_n,
\en
where we use a shorter notation for the multinomial coefficients
$$
\binom{m}{\gamma} = \frac{m!}{m_1! \dots m_d!}, \quad \gamma = (m_1,\dots,m_d), \ \
m_1 + \dots + m_d = m,
$$
and perform the summation over all multi-indices $\gamma_i$ such that
$\gamma_1 + \dots + \gamma_n = \gamma$. The latter implies that
$|\gamma_1| + \dots + |\gamma_n| =  |\gamma| = m$. 

On the one hand, the $\gamma$-derivative of the function $f(x) = (x_1 \dots x_d)^{4p}$ 
at the point $x = e = (1,\dots,1)$ is given by
$$
D^\gamma f(e) \, = \, \prod_{i=1}^d 4p (4p-1) \dots (4p-m_i+1) \, \leq \,
\prod_{i=1}^d (4p)^{m_i} \, = \, (4p)^m.
$$

On the other hand, since $f(x)$ may be written as a product of $n = 2p$ equal functions 
$$
f_i(x) = g(x) = x_1^2 \dots x_d^2,
$$ 
(4.4) is applicable, and we get
\be
\sum \binom{|\gamma_1|}{\gamma_1} \dots  \binom{|\gamma_{2p}|}{\gamma_{2p}}\, 
D^{\gamma_1} g(e) \dots D^{\gamma_{2p}} g(e) \, \leq \, (4p)^m
\en
with summation as before. It should be clear that $D^\gamma g(e) \geq 0$ 
for any multi-index $\gamma$.

Let us now apply (4.4) to the function 
\be
f(x) = (y (1 - y))^p, \quad y = y(x) =  \frac{4\,|x|^2 - 1}{3},
\en
which is also product of $2p$ functions 
$$
f_i(x) = y \ \ {\rm for} \ 1 \leq i \leq p, \ \ \ {\rm and} \ \ \ 
f_i(x) = 1-y \ \ {\rm for} \ p+1 \leq i \leq 2p. 
$$
Given a multi-index $\gamma = (m_1,\dots,m_d)$, note that $D^\gamma |x|^2 = 0$ 
if $m_i \geq 1$ for at least two distinct values of $i$. In the other case
$\gamma = me_i$ where $e_i$ is the $i$-th vector in the canonical basis in $\R^d$, 
\[
D^{me_i} |x|^2 = \left\{
\begin{array}{cc}
\!\!|x|^2, & \ {\rm if} \ m = 0, \\
\!\!2x_i,   & \ {\rm if} \ m = 1, \\
\!\!2,      & \ {\rm if} \ m = 2, \\
\!\!0,      & \ {\rm if} \ m \geq 3,
\end{array}
\right. \quad
D^{me_i} g(e) = \left\{
\begin{array}{cc}
\!\!1, & \ {\rm if} \ m = 0, \\
\!\!2,   & \ {\rm if} \ m = 1, \\
\!\!2,      & \ {\rm if} \ m =2, \\
\!\!0,      & \ {\rm if} \ m \geq 3.
\end{array}
\right.
\]
Hence $\big|D^\gamma |x|^2\big| \leq D^\gamma g(e)$ on the unit ball $|x| \leq 1$
for any multi-index $\gamma$. Since also $D^\gamma |x|^2 = 0$ whenever $|\gamma| \geq 3$,
we get
$$
|D^\gamma y| \, = \, \Big(\frac{4}{3}\Big)^{|\gamma|}\, \big|D^\gamma |x|^2\big|
 \, \leq \, \Big(\frac{4}{3}\Big)^2 D^\gamma g(e), \quad |\gamma|>0.
$$
A similar bound holds true for $D^\gamma (1-y)$ and also in the case $|\gamma| = 0$ for 
$\frac{1}{2} \leq |x| \leq 1$, since then $0 \leq y \leq 1$, while $D^\gamma g(e) = g(e) = 1$. 
As a consequence of (4.4)-(4.5), the function $f(x)$ in (4.6) satisfies in the region $\frac{1}{2} \leq |x| \leq 1$
\begin{eqnarray}
\big|D^\gamma f(x)\big| 
 & \leq &
\sum \binom{|\gamma_1|}{\gamma_1} \dots  \binom{|\gamma_{2p}|}{\gamma_{2p}}\, 
|D^{\gamma_1} y| \dots |D^{\gamma_p} y|\, |D^{\gamma_{p+1}} (1-y)| \dots
|D^{\gamma_{2p}} (1-y)| \nonumber \\
 & \leq &
\Big(\frac{4}{3}\Big)^{4p} \, \sum \binom{|\gamma_1|}{\gamma_1} \dots  \binom{|\gamma_{2p}|}{\gamma_{2p}}\, 
D^{\gamma_1} g(e) \dots D^{\gamma_{2p}} g(e) \nonumber \\
 & \leq &
\Big(\frac{4}{3}\Big)^{4p}\, (4p)^m.
\end{eqnarray}

Let us now return to the function $\psi(x) = 1 - v(y)$ as defined in (4.2)-(4.3)
and assume for definiteness that $m_1 \geq 1$. Since
$$
\partial_{x_1} \psi(x) = - \frac{8}{3A_p} x_1 f(x), \quad f(x) = (y (1 - y))^p,
$$
we have
$$
D^\gamma \psi(x) = - \frac{8}{3A_p} \big(x_1 D^\gamma f(x) + D^{\gamma - e_1} f(x)\big).
$$
Hence, by (4.7), we get for the region $\frac{1}{2} \leq |x| \leq 1$, 
$$
|D^\gamma \psi(x)| \leq \frac{16}{3A_p} \Big(\frac{4}{3}\Big)^{4p}\, (4p)^m, \quad 
0 \leq |\gamma| \leq m.
$$
We finally note that
$$
A_p = \int_0^1 (s(1-s))^p\,ds = \frac{p!^2}{(2p+1)!} \geq \frac{1}{4^p (2p+1)}.
$$
Thus, 
$$
\sup_{|x| \leq 1} |D^\gamma \psi(x)| \, \leq \,
4^p (2p+1) \frac{16}{3} \Big(\frac{4}{3}\Big)^{4p}\, (4p)^m.
$$
\qed

\vskip7mm
\section{{\bf Reduction to Compactly Supported Functions}}
\setcounter{equation}{0}

\vskip2mm
\noindent
Let $U_r$ ($r>0$) denote the collection of all $C^\infty$-smooth functions $u:\R^d \rightarrow \R$ 
which are supported on the ball $B_r = \{x \in \R^d: |x| \leq r\}$ and satisfy
the conditions (2.3)-(2.4) in Lemma 2.1. As a next step, we show that the supremum in definition (2.1) 
of the Zolotarev semi-norm $\|\lambda\|_{\zeta_p}$ may be restricted to the set $U_r$ at the 
expense of a small error for large values of the parameter $r$ under a $q$-th moment assumption
\be
\int_{\R^d} |x|^q\,d|\lambda|(x) \leq b^q \quad (b \geq 0, \ q>p).
\en

Define
\be
\|\lambda\|_{\zeta_p(r)} = \sup_{u \in U_r} \Big|\int_{B_r} u\,d\lambda\Big|.
\en 
As we noted in (2.6), any function $u$ on $\R^d$ satisfying (2.3)-(2.4) admits a polynomial
pointwise bound
\be
|u(x)| \leq \frac{1}{p!}\,d^{\frac{p}{2}}\,|x|^p, \quad x \in \R^d.
\en
Hence, for any signed measure $\lambda$ on $\R^d$,
$$
\|\lambda\|_{\zeta_p(r)} \, \leq \, \frac{r^p}{p!}\,d^{\frac{p}{2}}\,\|\lambda\|_{\rm TV}.
$$
Thus, for the finiteness of this semi-norm, the moment condition (2.2) is not needed.
However, it is needed to bound the Zolotarev norm $\|\lambda\|_{\zeta_p}$.

\vskip5mm
{\bf Lemma 5.1.} {\sl Let the signed measure $\lambda$ have zero mixed moments up to order $p-1$
and satisfy $(5.1)$ with $q > p$ and $b \geq 0$. Then with some absolute constant $c>0$
\be
\|\lambda\|_{\zeta_p}  \leq c^p d^{p/2}
\, \|\lambda\|_{\zeta_p(r)} + \frac{1}{p!}\,d^{p/2}\, \Big(\frac{2}{r}\Big)^{q-p}\,b^q.
\en
}

\vskip3mm
{\bf Proof.} The assumption on $\lambda$ includes the requirement that it has total mass zero.

Following Lemma 2.1, take a smooth function $u$ on $\R^d$ satisfying (2.3)-(2.4)
and define a new function
\be
u_r(x) = u(x) \psi(x/r), \quad x \in \R^d,
\en
where the $C^p$-smooth function $\psi:\R^d \rightarrow [0,1]$ is described in Lemma 4.1.
From the properties 1)-2) in this lemma it follows that 
\be
|u_r(x) - u(x)| \leq |u(x)| \quad {\rm for} \ \frac{r}{2} \leq |x| \leq r,
\en
\be
u_r(x)- u(x) = 0 \ {\rm for} \ |x| \leq \frac{r}{2}, \quad
u_r(x) = 0 \ {\rm for} \ |x| \geq r,
\en
\be
D^\alpha u_r(0) = 0 \ \ {\rm whenever} \ |\alpha| \leq p-1.
\en
Using property 3)  of this lemma, one can also show that
\be
|D^\alpha u_r(x)| \leq C \ \ {\rm whenever} \ |\alpha| = p
\en
for some constants $C = C_{p,d}$ depending on $(p,d)$.
This dependence will be clarified later on.

Using the moment assumption (5.1) together with (5.3) and (5.6)-(5.7), we have
\bee
\Big|\int_{\R^d} u_r\,d\lambda - \int_{\R^d} u\,d\lambda \Big| 
 & \leq &
\int_{|x| > r/2} |u(x)|\,d|\lambda|(x) \\
 & \leq &
\frac{1}{p!}\,d^{p/2} \int_{|x| > r/2} |x|^p\,d|\lambda|(x) \\
 & \leq &
\frac{1}{p!}\,d^{p/2} \int_{|x| > r/2} \frac{|x|^q}{(r/2)^{q-p}}\,d|\lambda|(x) \, \leq \,
\frac{1}{p!}\,d^{q/2}\, \Big(\frac{2}{r}\Big)^{q-p}\,b^q.
\ene
Thus,
\bee
\Big|\int_{\R^d} u\,d\lambda \Big| 
 & \leq &
\Big|\int_{\R^d} u_r\,d\lambda\Big| + \frac{1}{p!}\,d^{q/2}\, \Big(\frac{2}{r}\Big)^{q-p}\,b^q \\
 & \leq &
C\, \|\lambda\|_{\zeta_p(r)} + \frac{1}{p!}\,d^{q/2}\, \Big(\frac{2}{r}\Big)^{q-p}\,b^q, 
\ene
where we used the property that, the function $\frac{1}{C}\,u_r$ participates in the 
supremum (5.2), according to (5.8)-(5.9). Taking the supremum on the left-hand side 
over all admissible functions $u$, we get the required inequality (5.4) with a constant $C$ 
instead of the explicit constant $c^p d^{p/2}$.

Let us now estimate the constant $C$.
By the Newton binomial differentiation formula applied to the product in (5.5),
for any multi-index $\alpha = (p_1,\ldots,p_d)$ of length $|\alpha| = p$, we obtain
$$
D^\alpha u_r(x) \, = \,
\sum \binom{p_1}{k_1}  \cdots \binom {p_d}{k_d}  \, 
D^\beta u(x) \, r^{-(p - |\beta|)}\,D^{\alpha - \beta} \psi(x/r),
$$
where the summation is performed over all multi-indices $\beta = (k_1,\dots,k_d)$ such that 
$k_i \leq p_i$ for all $i \leq d$, i.e. $\beta \leq \alpha$ using the componentwise comparison.
By Lemma 2.2,
$$
|D^\beta u(x)| \leq \frac{1}{(p - |\beta|)!}\,d^{\frac{p-|\beta|}{2}}\,r^{p-|\beta|}, \quad
|x| \leq r.
$$
Hence we get a uniform bound
\be
|D^\alpha u_r(x)| \, \leq  \,
\sum_{\beta \leq \alpha} \binom{p_1}{k_1}  \cdots \binom {p_d}{k_d}  \, 
\frac{1}{(p - |\beta|)!}\,d^{\frac{p-|\beta|}{2}}\, \sup_{|x| < 1} |D^{\alpha - \beta} \psi(x)|,
\en
which does not depend on $r$.

Now, let us recall that, by Lemma 4.1, with some absolute constant $c_0 > 0$
$$
\sup_{|x| \leq 1} |D^\gamma \psi(x)| \leq c_0^p\,p^{|\gamma|}, \quad 0 \leq |\gamma| \leq p.
$$
Applying this bound in (5.10) with $\gamma = \alpha - \beta$ and putting $m = |\gamma| = p - |\beta|$, 
we get
$$
|D^\alpha u_r(x)| \, \leq  \, c_0^p \,d^{p/2}\, 
\sum_{\gamma \leq \alpha} \binom{p_1}{k_1}  \cdots \binom {p_d}{k_d}  \, \frac{1}{m!}\,p^m.
$$
If we fix $m = 0,1,\dots,p$ and extend the sum to all $\beta = (k_1,\dots,k_d)$ such that 
$k_i \leq p_i$ for all $i \leq d$, we will get the value $2^{p_1} \dots 2^{p_d} = 2^p$. Hence
$$
|D^\alpha u_r(x)| \, \leq  \, (2c_0)^p \,d^{p/2}\, \sum_{m=0}^p \frac{1}{m!}\,p^m <
(2c_0)^p \,d^{p/2}\,e^p.
$$
Thus $C \leq (2c_0 e)^p \,d^{p/2}$.
\qed

\vskip7mm
\section{{\bf Fourier Transforms}}
\setcounter{equation}{0}

\vskip2mm
\noindent
Any integrable compactly supported function $u$ on $\R^d$ has 
the Fourier transform
\be
\widehat u(t) = \int_{\R^d} e^{i \left<t,x\right>} u(x)\,dx, \quad t \in \R^d,
\en
which represents a $C^\infty$-smooth function. 
Towards the proof of Theorem 1.2 let us state now the following integrability property.

\vskip5mm
{\bf Lemma 6.1.} {\sl For any $u:\R^d \rightarrow \R$ which is supported on the ball $B_r$ 
and has partial derivatives of order $p$ bounded by $1$ in absolute value,
\be
\int_{\R^d} |\widehat u(t)|^2\, |t|^{2p}\,dt \, \leq \, d^p\omega_d\, (2\pi r)^d.
\en
}

\vskip2mm
In particular, this inequality holds true for any function $u$ in $U_r$.

\vskip5mm
{\bf Proof.} First assume that $\widehat u(t) = O(1/|t|^n)$ as $|t| \rightarrow \infty$ for any 
integer $n \geq 1$. Then (6.1) may be inverted via the inverse Fourier transform
$$
u(x) = \frac{1}{(2\pi)^d} \int_{\R^d} e^{-i \left<t,x\right>}\, \widehat u(t)\,dt.
$$
In particular, $u$ is $C^\infty$-smooth on $\R^d$. Moreover, this equality may be
differentiated $p$ times along every coordinate $x_k$ to represent the corresponding partial
derivatives as
$$
\partial_{x_k}^p u(x) =
\frac{(-i)^p}{(2\pi)^d} \int_{\R^d} e^{-i \left<t,x\right>}\,t_k^p\, \widehat u(t)\,dt, \quad k = 1,\dots,d,
$$
where $t = (t_1,\dots,t_d)$. Hence, by the Plancherel theorem,
$$
\int_{\R^d} t_k^{2p}\, |\widehat u(t)|^2\,dt \, = \,
(2\pi)^d \int_{\R^d} \big(\partial_{x_k}^p u(x)\big)^2\,dx \, = \,
(2\pi)^d \int_{B_r} \big(\partial_{x_k}^p u(x)\big)^2\,dx,
$$
where we also used the assumption that $u$ and hence its partial derivatives are vanishing 
outside the ball $B_r$. The last integrand is bounded by 1, and
summing over all $k \leq d$, we get
$$
\int_{\R^d} (t_1^{2p} + \dots + t_d^{2p})\,|\widehat u(t)|^2\, dt \, \leq \,
d (2\pi)^d \cdot \omega_d r^d.
$$
Using $t_1^{2p} + \dots + t_d^{2p} \geq d^{-(p-1)}\,|t|^{2p}$, we arrive at the desired inequality (6.2).

In the general case, a smoothing argument can be used. By the well-known theorem
of Ingham \cite{I}, there exists a probability density $w$ 
on $\R^d$ which is supported on $B_1$ and has characteristic function 
$\widehat w(t)$ satisfying $\widehat w(t) = O(1/|t|^n)$ as $|t| \rightarrow \infty$,
for any fixed $n>0$.
Given $\ep>0$, the probability density $w_\ep(x) = \ep^{-d} w(x/\ep)$ is supported 
on the ball $B_\ep$ and has characteristic function $\widehat w_\ep(t) = \widehat w(\ep t)$. 
Consider the convolution
$$
u_\ep(x) = (u * w_\ep)(x) = \int_{\R^d} u(x-y)\,w_\ep(y)\,dy, \quad x \in \R^d.
$$
This function is $C^\infty$-smooth and is supported on $B_{r+\ep}$. Moreover,
the above equality may be differentiated $p$ times along the $k$-th coordinate to yield
$$
\partial_{x_k}^p u(x) = \big(\partial_{x_k}^p u * w_\ep\big)(x) = 
\int_{\R^d} \partial_{x_k}^p u(x-y)\,w_\ep(y)\,dy,
$$
implying that $|\partial_{x_k}^p u(x)| \leq 1$ for all $x \in \R^d$.

Since $u$ is bounded and compactly supported, it follows from (6.1) that
$\sup_{t \in \R^d} |\widehat u(t)| < \infty$. Hence, the Fourier transform of $u_\ep$ satisfies
$\widehat u_\ep(t) = \widehat u(t) \widehat w(\ep t) = O(1/|t|^n)$ as $|t| \rightarrow \infty$.
Thus, one may apply the previous step to the function $u_\ep$ which gives
$$
\int_{\R^d} |\widehat u(t)|^2\, |\widehat w(\ep t)|^2\,|t|^{2p}\,dt \, \leq \, 
d^p\omega_d\, (2\pi\, (r+\ep))^d.
$$
It remains to send $\ep \rightarrow 0$ and apply Fatou's lemma together with
$\widehat w(\ep t) \rightarrow 1$ as $\ep \rightarrow 0$.
\qed

\vskip4mm
Next, let us connect the Kantorovich distance with Fourier-Stieltjes transforms
(or multivariate characteristic functions)
$$
f(t) = \int_{\R^d} e^{i \left<t,x\right>} d\mu(x), \quad
g(t) = \int_{\R^d} e^{i \left<t,x\right>} d\nu(x) \qquad (t \in \R^d).
$$

\vskip5mm
{\bf Lemma 6.2.} {\sl Given two probability measures $\mu$ and $\nu$ in ${\mathfrak P}_p$
with characteristic functions $f$ and $g$, for any $t \in \R^d$,
\be
|f(t) - g(t)| \leq 2|t|^p\,\sup_{|\theta| = 1}  \zeta_p(\mu_\theta,\nu_\theta).
\en
}
\vskip2mm
{\bf Proof.} In dimension one, using the property that the function $u_t(x) = t^{-p}\,\cos(itx)$
($t \neq 0$) has a $p$-th derivative bounded by 1 in absolute value,
it follows from (1.5) that
$$
|{\rm Re}(f(t)) - {\rm Re}(g(t))| \leq |t|^p\,\zeta_p(\mu,\nu).
$$
By a similar argument,
$$
|{\rm Im}(f(t)) - {\rm Im}(g(t))| \leq |t|^p\,\zeta_p(\mu,\nu),
$$
so that
$$
|f(t) - g(t)| \leq 2|t|^p\,\zeta_p(\mu,\nu).
$$
In dimension $d$, note that, for any $\theta \in \R^d$, the functions
$r \rightarrow f(r\theta)$ and $r \rightarrow g(r\theta)$
represent the characteristic functions for the images $\mu_\theta$ and $\nu_\theta$ 
of $\mu$ and $\nu$ under the linear map $x \rightarrow \left<x,\theta\right>$. Hence
$$
|f(r\theta) - g(r\theta)| \leq 2|r|^p\,\zeta_p(\mu_\theta,\nu_\theta).
$$
Taking the supremum of both sides over all $\theta \in S^{d-1}$ with $r = |t|$,
we arrive at (6.3).
\qed

\vskip7mm
\section{{\bf Proof of Theorem 1.2}}
\setcounter{equation}{0}

\noindent
Let $X$ and $Y$ be random vectors distributed according to $\mu$ and $\nu$ with finite
moments
$$
\|X\|_q = (\E\,|X|^q)^{1/q}, \quad \|Y\|_q = (\E\,|Y|^q)^{1/q}.
$$
With $b = \max(\|X\|_q,\|Y\|_q)$ the inequality (1.4) is homogeneous with respect to $(X,Y)$, 
so one may assume that $b=1$. In particular, $\|X\|_p \leq 1$ and $\|Y\|_p \leq 1$.
As a consequence, applying Corollary 2.4 with the signed measure
$$
\lambda = \mu - \nu
$$ 
and using $|\lambda| \leq \mu + \nu$, we have
\be
\zeta_p(\mu,\nu) \leq \frac{d^{\frac{p}{2}}}{p!} \big(\|X\|_p^p + \|Y\|_p^p\big) \leq 
\frac{2 d^{\frac{p}{2}}}{p!}.
\en

Recall that in (3.2) we introduced the convolutions 
$$
\lambda_\ep = \lambda * \kappa_\ep  = \mu * \kappa_\ep - \nu * \kappa_\ep
$$ 
with parameter $\ep>0$. Since $ \|\lambda\|_{\rm TV} \leq 2$, an application
of Corollary 3.3 yields
$$
\|\lambda_\ep - \lambda\|_{\zeta_p} \, \leq \, \frac{(8d\ep)^p}{(p-1)!}, 
$$
so that, by the triangle inequality for the Zolotarev semi-norm,
\be 
\|\lambda\|_{\zeta_p} \, \leq \,\|\lambda_\ep\|_{\zeta_p} +  \frac{(8d\ep)^p}{(p-1)!}. 
\en

Let us now focus on the convolved measures $\lambda_\ep$ and estimate their 
Zolotarev semi-norms. As was mentioned before, every $\lambda_\ep$ has zero mixed moments up 
to order $p-1$. An application of Lemma 5.1 to $\lambda_\ep$ in place of $\lambda$ gives
\be
\|\lambda_\ep\|_{\zeta_p}  \leq 
c^p d^{p/2} \, \|\lambda_\ep\|_{\zeta_p(r)} + \frac{1}{p!}\,d^{p/2}\, \Big(\frac{2}{r}\Big)^{q-p}\,b_\ep^q
\en
with some absolute constant $c$, where
$$
b_\ep = \Big(\int_{\R^d} |x|^q\,d|\lambda_\ep|(x)\Big)^{1/q}.
$$
In order to estimate this quantity, we use for the total variation measure the bound
$$
|\lambda_\ep| \leq \mu * |\kappa_\ep| + \nu * |\kappa_\ep|.
$$
Note that $\|\kappa_\ep\|_{\rm TV} = \|\kappa\|_{\rm TV}$.
Let $Z_\ep$ denote a random vector with distribution $\frac{1}{\|\kappa\|_{\rm TV}}\,|\kappa_\ep|$,
assuming that it is independent of $X$. Since $\kappa_\ep$ and therefore $|\kappa_\ep|$ are 
supported on the ball $B_\ep$, we have $\|Z_\ep\|_q \leq \ep$. Hence, by the Minkowski inequality,
\bee
\int_{\R^d} |x|^q\,d\,(\mu * |\kappa_\ep|)(x)
 & = &
\int_{\R^d} \int_{\R^d} |x-y|^q\,d\mu(x)\, d|\kappa_\ep|(y) \\
 & = & 
\|\kappa_\ep\|_{\rm TV}\,\E\,|X + Z_\ep|^q \\
 & = &
 \|\kappa\|_{\rm TV}\, \|X + Z_\ep\|_q^q \\
 & \leq &
\|\kappa\|_{\rm TV}\, \big(\|X\|_q + \|Z_\ep\|_q\big)^q \, \leq \, \|\kappa\|_{\rm TV}\, (1 + \ep)^q.
\ene
The total variation norm of $\kappa$ was estimated in Lemma 3.2, which gives
$$
\int_{\R^d} |x|^q\,d\,(\mu * |\kappa_\ep|)(x) \leq p\,2^{3p-1}\,(1 + \ep)^q.
$$
A similar bound holds for $\nu$, and therefore $b_\ep^q \leq  p\,2^{3p}\,(1 + \ep)^q$.
Applying this in (7.3), we get
$$
\|\lambda_\ep\|_{\zeta_p}  \leq c^p d^{p/2}  \, \|\lambda_\ep\|_{\zeta_p(r)} + 
\frac{8^p}{(p-1)!}\,d^{p/2}\, \Big(\frac{2}{r}\Big)^{q-p}\,(1+\ep)^q.
$$
Thus, by (7.2),
\begin{eqnarray}
\|\lambda\|_{\zeta_p} 
 & \leq & 
c^p d^{p/2} \, \|\lambda_\ep\|_{\zeta_p(r)} \nonumber \\
 & & \hskip-5mm + \
\frac{8^p}{(p-1)!}\,d^{p/2}\, \Big(\frac{2}{r}\Big)^{q-p}\,(1+\ep)^q + \frac{(8d\ep)^p}{(p-1)!}.
\end{eqnarray}

Next, we need to estimate the semi-norm $\|\lambda_\ep\|_{\zeta_p(r)}$. To this aim, introduce 
the Fourier transform $h(t)$ of the measure $\kappa$, whose density was denoted by $w$ 
in Lemma 3.2. The convolved measures $\mu_\ep$ and $\nu_\ep$ have continuous densities 
$p_\ep$ and $q_\ep$ with Fourier transforms
$$
f_\ep(t) = f(t) h(\ep t), \quad g_\ep(t) = g(t) h(\ep t) \qquad (t \in \R^d),
$$
where $f$ and $g$ denote the characteristic functions of $X$ and $Y$.
Hence, given a function $u$ in $U_r$, one may write
$$
\int_{\R^d} u\,d\lambda_\ep = \int_{\R^d} u(x)\,(p_\ep(x) - q_\ep(x))\,dx.
$$
We are in position to apply the Plancherel theorem and 
rewrite the last integral as
$$
\frac{1}{(2\pi)^d} \int_{\R^d} \widehat u(t)\,(\bar f(t) - \bar g(t))\,h(\ep t)\,dt
$$
(due to the symmetry of $\kappa$, the function $h(t)$ is real-valued).
Thanks to Lemma 6.2,
\be
\Big|\int_{\R^d} u\,d\lambda_\ep\Big| \leq 
\frac{2M}{(2\pi)^d} \int_{\R^d} |\widehat u(t)|\,|t|^p\,|h(\ep t)|\,dt,
\en
where 
$$
M = \sup_{|\theta| = 1} \zeta_p(\mu_\theta,\nu_\theta).
$$ 

To proceed, apply the Plancherel theorem once more together with the uniform
bound (3.5) on $w(x)$ from Lemma 3.2 to get that
\bee
\int_{\R^d} |h(\ep t)|^2\,dt 
 & = &
\ep^{-d} \int_{\R^d} |h(t)|^2\,dt \\
 & = &
(2\pi)^d\, \ep^{-d} \int_{\R^d} w(x)^2\,dx \, \leq \, 
(2\pi)^d\, \ep^{-d} \cdot \frac{1}{\omega_d}\,p^{d+1}\, 8^p.
\ene
Applying Cauchy's inequality in (7.5) and recalling Lemma 6.1, we then obtain that
\bee
\Big|\int_{\R^d} u\,d\lambda_\ep\Big|  
 & \leq &
\frac{2M}{(2\pi)^d} \bigg(\int_{\R^d} |\widehat u(t)|^2\,|t|^{2p}\,dt\bigg)^{1/2}
\bigg(\int_{\R^d} |h(\ep t)|^2\,dt\bigg)^{1/2} \\
 & \leq &
\frac{2M}{(2\pi)^d}\,\Big(d^p \omega_d\,(2\pi r)^d\Big)^{1/2} 
\Big((2\pi)^d\, \ep^{-d} \cdot \frac{1}{\omega_d}\,p^{d+1}\, 8^p\Big)^{1/2} \\
 & = &
2M\, (8d)^{\frac{p}{2}}\, \Big(\frac{pr}{\ep}\Big)^{\frac{d}{2}} \sqrt{p}.
\ene
Taking the supremum over all $u \in U(r)$ on the left-hand side leads
to the similar bound for $\|\lambda_\ep\|_{\zeta_p(r)}$, and using this in (7.4) we are led to
$$
\|\lambda\|_{\zeta_p} \, \leq \,c^p d^{p/2}  \,(8d)^{\frac{p}{2}}\,M\, 
\Big(\frac{pr}{\ep}\Big)^{\frac{d}{2}} \sqrt{p} + 
\frac{8^p}{(p-1)!}\,d^{p/2}\, \Big(\frac{2}{r}\Big)^{q-p}\,(1+\ep)^q +
\frac{(8d\ep)^p}{(p-1)!}.
$$
Here the terms $8^{p/2}$ and $\sqrt{p}$ can be absorbed in $c^p$, so that the above simplifies to
$$
\|\lambda\|_{\zeta_p} \, \leq \,(cd)^p\,M\, \Big(\frac{pr}{\ep}\Big)^{\frac{d}{2}} + 
\frac{8^p}{(p-1)!}\,d^{p/2}\, \Big(\frac{2}{r}\Big)^{q-p}\,(1+\ep)^q + \frac{(8d\ep)^p}{(p-1)!}. 
$$

To simplify optimization over free variables $r>0$ and $\ep>0$, let us require that that
$\ep \leq 1$. Then for some absolute constant $c>0$ we have
$$
c^p\, \|\lambda\|_{\zeta_p} \, \leq \, d^p M\,\Big(\frac{pr}{\ep}\Big)^{\frac{d}{2}} + 
\frac{1}{(p-1)!}\,d^{p/2}\, \Big(\frac{4}{r}\Big)^{q-p} + \frac{(d\ep)^p}{(p-1)!}.
$$
Let us then replace $r$ with $4s$ in the above inequality to get 
$$
c^p\,\|\lambda\|_{\zeta_p} \, \leq \, d^p M\, \Big(\frac{4p s}{\ep}\Big)^{\frac{d}{2}} + 
\frac{1}{(p-1)!}\,d^{p/2}\, \Big(\frac{1}{s}\Big)^{q-p} + \frac{(d\ep)^p}{(p-1)!}.
$$
Further replacement of $\ep$ with $4p\delta$ and choosing a smaller value of $c>0$ 
leads to
$$
c^p\,
\|\lambda\|_{\zeta_p} \, \leq \, d^p M\, \Big(\frac{s}{\delta}\Big)^{\frac{d}{2}} + 
\frac{1}{p!}\,d^{p/2}\, \Big(\frac{1}{s}\Big)^{q-p} + \frac{(d p)^p}{p!} \delta^p.
$$
By Stirling's formula, the latter fraction does not exceed $4^p$, and choosing
a smaller constant $c$, we get
$$
c^p\,
\|\lambda\|_{\zeta_p} \, \leq \, d^p M\, \Big(\frac{s}{\delta}\Big)^{\frac{d}{2}} + 
\frac{1}{p!}\,d^{p/2}\,s^{p-q} + (d\delta)^p.
$$

Here, equalizing the terms $M\,(\frac{s}{\delta})^{\frac{d}{2}}$ and $s^{p-q}$,
we find the unique value 
$$
s = \bigg(\frac{\delta^{\frac{d}{2}} }{M}
\bigg)^{\frac{1}{ q-p + \frac{d}{2} }}
$$ 
for which the above yields
\be
c^p\,\|\lambda\|_{\zeta_p} \leq d^p( A \delta^{-\sigma} + \delta^p),
\en
where
$$
\sigma = \frac{(q-p)\,\frac{d}{2}}{q-p + \frac{d}{2}}, \quad
A = M^{\frac{q-p}{q-p + \frac{d}{2}}}.
$$
The choice $\delta = A^{\frac{1}{\sigma + p}}$ in (7.6) leads to
\be
c^p\,\|\lambda\|_{\zeta_p} \, \leq \, 2d^p A^{\frac{p}{\sigma + p}}.
\en
Recall that this holds provided that $4p\delta \leq 1$,
that is, if $A^{\frac{1}{\sigma + p}} \leq \frac{1}{4p}$. 

In the other case where $A^{\frac{1}{\sigma + p}} \geq \frac{1}{4p}$, 
the relation (7.7) would follow from
$$
c^p\,\|\lambda\|_{\zeta_p} \leq \frac{2d^p }{(4p)^p}. 
$$
And this is indeed true for a suitable constant $c$ in view of the general bound (7.1).
Thus, in both case we proved that
$$
\|\lambda\|_{\zeta_p} \, \leq \, (cd)^p\, A^{\frac{p}{\sigma + p}}
$$
with some absolute constant $c>0$. It remains to note that
$$
A^{\frac{p}{\sigma + p}} = M^\beta, \quad
\beta = \frac{2}{2 + \frac{dq}{p(q-p)}}.
$$
\qed


\end{document}